\documentclass[]{article}
\usepackage{amssymb}
\usepackage{amsmath}
\usepackage{amsthm}
\usepackage{hyperref}

\newtheorem{Theorem}{Theorem}

\newtheorem{Definition}{Definition}

\newcommand{\pr}{\mbox{pr}_1\,}

\title {Stabilizers  of functional Menger systems}
\author {Wieslaw A. Dudek and Valentin S. Trokhimenko}
\date {}

\begin{document}
\sloppy \maketitle

\begin{abstract}
A functional Menger system is a set of $n$-place functions
containing $n$ projections and closed under the so-called Menger's
composition of $n$-place functions. We give the abstract
characterization for subsets of these functional systems which
contain functions  having one common fixed point.
\end{abstract}

\section{Introduction}
Investigation of multiplace functions by algebraic methods plays a
very important role in modern mathematics were we consider various
operations on sets of functions, which are naturally defined. The
basic operation for functions is superposition (composition), but
there are  some other naturally defined operations, which are also
worth of consideration. For example, the operation of
set-theoretic intersection and the operation of projections  (see
for example \cite{Dudtro, Dudtro1, Dudtro2, schtr, Tr-1}). The
central role in these study play sets of functions with fixed
points. The study of such sets for functions of one variable was
initiated by B. M. Schein in \cite{146} and \cite{SchLectures}.
Next, for sets of functions of $n$ variables, was continued by
V.S. Trokhimenko (see \cite{ Tr-3, Tr-4}).

In this paper, we consider the sets of $n$-place functions
containing $n$-projections and closed under the so-called Menger's
composition of $n$-place functions. For such functional systems we
give the abstract characterization for subsets of functions having
one common fixed point.

\section{Preliminaries}

Let $A^n$ be the $n$-th Cartesian product of a set $A$. Any
partial map from $A^n$ into $A$ is called an \textit{$n$-place
function} on $A$. The set of all such maps is denoted by
$\mathcal{F}(A^n,A)$. On ${\mathcal F}(A^n,A)$ we define one
$(n+1)$-ary superposition $\mathcal{O}\colon (f,g_1,\ldots,
g_n)\mapsto f[g_1\ldots g_n]$, called the {\it Menger's
composition}, and $n$ unary operations $\mathcal{R}_i\colon
f\mapsto\mathcal{R}_if$, $i\in\overline{1,n}=\{1,\ldots,n\}$
putting
\begin{eqnarray}
f [g_1\ldots g_n](a_1,\ldots,a_n)=f(g_1(a_1,\ldots,a_n),\ldots,
g_n(a_1,\ldots,a_n)), \label{super} \\[4pt]
\mathcal{R}_if(a_1,\ldots,a_n)=a_i, \ \mbox{where} \ \
(a_1,\ldots,a_n) \in\pr f
\end{eqnarray}
for $f,g_1,\ldots,g_n\in\mathcal{F}(A^n,A)$, $(a_1,\ldots,a_n)\in
A^n$, where $\pr f$ denotes the domain of a function $f$. It is
assumed the left and right hand side of equality (\ref{super}) are
defined, or not defined, simultaneously. Algebras of the form
$(\Phi,\mathcal{O},\mathcal{R}_1,\ldots,\mathcal{R}_n)$, where
$\Phi\subset\mathcal{F}(A^n, A)$, are called \textit{functional
Menger systems of $n$-place functions}. Algebras of the form
$(\Phi,\mathcal{O},\cap,\mathcal{R}_1,\ldots,\mathcal{R}_n)$,
where $\cap$ is a set-theoretic intersection, are called
\textit{functional Menger $\cap$-algebras of $n$-place functions}.
In the literature such algebras are also called \textit{functional
Menger $\cal{P}$-algebras} (see \cite{Dudtro} and \cite{Tr-4}).

Let $a$ be some fixed element (point) of $A$. The
\textit{stabilizer of $a$} is the set $H^{a}_{\Phi}$ of such
functions from $\Phi$ for which $a$ is a fixed point, i.e., the
set
\[
H^{a}_{\Phi}=\{f\in\Phi\,|\,f(a,\ldots,a)=a\}.
 \]

Let $(G,o)$ be a nonempty set with one $(n+1)$-ary operation
\[o\colon\ (x_0,x_1,\ldots,x_n)\mapsto x_0[x_1\ldots x_n].\]
 An abstract algebra
$\mathcal{G}=(G,o, R_1,\ldots,R_n)$ of type $(n+1,1,\ldots,1)$ for
all $i,k\in\overline{1,n}$ satisfying the following axioms:
\[\begin{array}{cl}
\mathbf{A_1}\colon& \quad x[y_1\ldots y_n][z_1\ldots z_n]=
x[y_1[z_1\ldots z_n]\ldots y_n[z_1\ldots z_n]], \\[4pt]
\mathbf{A_2}\colon& \quad x[R_1x\ldots R_nx]=x, \\[4pt]
\mathbf{A_3}\colon& \quad x[\bar{u}\,|_iz][R_1y\ldots R_ny]=
x[\bar{u}\,|_iz[R_1y\ldots R_ny]], \\[4pt]
\mathbf{A_4}\colon& \quad R_ix[R_1y\ldots R_ny]= (R_ix)[R_1y\ldots
R_ny], \\[4pt]
\mathbf{A_5}\colon& \quad x[R_1y\ldots R_ny][R_1z\ldots R_nz]=
x[R_1z\ldots R_nz][R_1y\ldots R_ny], \\[4pt]
\mathbf{A_6}\colon& \quad R_ix[y_1\ldots y_n]=R_i(R_kx)[y_1\ldots y_n], \\[4pt]
\mathbf{A_7}\colon& \quad (R_ix)[y_1\ldots y_n]=
y_i[R_1(x[y_1\ldots y_n])\ldots R_n(x[y_1\ldots y_n])],
\end{array}
\]
where $x[\bar{u}\,|_iz]$ means $x[u_1\ldots
u_{i-1}z\,u_{i+1}\ldots u_n]$, is called a \textit{functional
Menger system of rank $n$}.

An algebra $\mathcal{G}_{\curlywedge}=(G,o,\curlywedge,
R_1,\ldots,R_n)$ of type $(n+1,2,1,\ldots,1)$, where
$(G,o,R_1,\ldots,R_n)$ is a functional Menger system of rank $n$
and $(G,\curlywedge)$ is a semilattice, is called a
\textit{functional Menger $\curlywedge$-algebra of rank $n$} if it
satisfies the identities:
\[
\begin{array}{cl}
\mathbf{A_8}\colon& \quad x\curlywedge y[R_1z\ldots R_nz]=
(x\curlywedge y)[R_1z\ldots R_nz], \\[4pt]
\mathbf{A_9}\colon& \quad x\curlywedge y=x[R_1(x\curlywedge
y)\ldots R_n(x\curlywedge y)], \\[4pt]
\mathbf{A_{10}}\colon& \quad (x\curlywedge y)[z_1\ldots z_n]=
x[z_1\ldots z_n]\curlywedge y[z_1\ldots z_n].
\end{array}
 \]

Any Menger algebra of rank $n$, i.e., an abstract groupoid $(G,o)$
with an $(n+1)$-ary operation satisfying $\mathbf{A_1}$ is
isomorphic to some set of $n$-place functions closed under
Menger's composition \cite{schtr}. Functional Menger
$\curlywedge$-algebras and Menger systems of rank $n$ are
isomorphic, respectively, to some functional Menger
$\cap$-algebras and Menger systems of $n$-place functions (see
\cite{Dudtro} and \cite{Tr-1}). Each homomorphism of such abstract
algebras into corresponding algebras of $n$-place functions is
called a \textit{representation by $n$-place functions}.
Representations which are isomorphisms are called
\textit{faithful}.

Let $(P_i)_{i\in I}$ be the family of representations of a Menger
algebra $(G,o)$ of rank $n$ by $n$-place functions defined on sets
$(A_i)_{i\in I}$, respectively. By the \textit{union} of this
family we mean the map $P \colon \, g \mapsto P(g)$, where $g\in
G$, and $P(g)$ is an $n$-place function on $A= \bigcup
\limits_{i\in I}A_i$ defined by
 \[
P(g)=\bigcup\limits_{i\in I}P_i(g).
 \]
If $A_i\cap A_j=\emptyset$ for all $i,j\in I$, $i\neq j$, then $P$
is called the \textit{sum} of $(P_i)_{i\in I}$ and is denoted by
$P=\sum\limits_{i\in I}P_i$. It is not difficult to see that the
sum of representations is a representation, but the union of
representations may not be a representation (see for example [1
$-$ 7]).

For every representation $P\colon G\rightarrow\mathcal{F}(A^n,A)$
of an algebra $(G,o)$ and every element $a\in G$ by $H^a_P$ we
denote the set of elements of $G$ corresponding to these $n$-place
functions for which $a$ is a fixed point, i.e.,
 \[
H^a_P=\{g\in G\,|\,P(g)(a,\ldots,a)=a\}.
 \]

Let $\mathcal{G}$ be a functional Menger system of rank $n$, $x$
-- an individual variable. By $T_n(G)$ we denote the set of
transformations $t\colon\,x\mapsto t(x)$ on $G$ such that:
\begin{enumerate}
\item[$(a)$] $x\in T_n(G)$,
\item[$(b)$] if $t(x)\in T_n(G)$, then $a[\,\bar{b}\,|_it(x)]\in
T_n(G)$ and $R_it(x)\in T_n(G)$ for all $a\in G$, $\bar{b}\in G^n$
and $i\in\overline{1,n}$,
\item[$(c)$] $T_n(G)$ contains only elements determined in $(a)$ and $(b)$.
\end{enumerate}

Let us remind that a nonempty subset $H$ of $G$ is called
\begin{itemize}
\item \textit{quasi-stable}, if for all $x\in G$
\[
x\in H\longrightarrow x[x\ldots x]\in H,
\]
\item \textit{$\curlywedge$-quasi-stable}, if for all $x\in G$
\[
x\in H\longrightarrow x[x\ldots x]\curlywedge x\in H,
\]
\item \textit{stable}, if for all $x,y_1,\ldots,y_n\in G$
\[
x,y_1,\ldots,y_n\in H\longrightarrow x[y_1\ldots y_n]\in H,
\]
\item \textit{$\curlywedge$-stable}, if for all $x,y\in G$
\[
x,y\in H\longrightarrow x\curlywedge y\in H,
\]
\item \textit{$l$-unitary}, if for every
$x,y\in G$
\[
x[y\ldots y]\in H\wedge\, y\in H\longrightarrow x\in H,
\]
\item \textit{$v$-unitary}, if for all $x,y_1,\ldots,y_n\in
G$
\[
x[y_1\ldots y_n]\in H\wedge\, y_1,\ldots,y_n\in H \longrightarrow
x\in H,
\]
\item a \textit{normal $v$-complex}, if for all $x,y\in
G$, $t\in T_n(G)$
\[
x,y\in H\wedge\, t(x)\in H\longrightarrow t(y)\in H,
\]
\item an \textit{$l$-ideal}, if for all $x,y_1,\ldots,y_n\in
G$
\[
(y_1,\ldots,y_n)\in G^{\,n}\setminus (G\setminus
H)^n\longrightarrow x[y_1\ldots y_n]\in H.
\]
\end{itemize}
A binary relation $\rho\subset G\times G$ is called
\begin{itemize}
\item \textit{stable}, if
\[
(x,y),(x_1,y_1),\ldots,(x_n,y_n)\in\rho\longrightarrow
(x[x_1\ldots x_n],y [y_1\ldots y_n])\in\rho
\]
for all $x,y,x_i,y_i\in G$, $i\in\overline{1,n}$,
\item \textit{$l$-regular}, if
\[
(x,y)\in\rho\longrightarrow (x[z_1\ldots z_n],y[z_1\ldots
z_n])\in\rho
\]
for all $x,y,z_i\in G$, $i\in\overline{1,n}$,
\item \textit{$v$-regular}, if
\[
(x_1,y_1),\ldots, (x_n,y_n)\in\rho\longrightarrow (z[x_1\ldots
x_n],z[y_1\ldots y_n])\in\rho
\]
for all $x_i,y_i,z\in G$, $i\in\overline{1,n}$,
\item \textit{$i$-regular}, where $i\in\overline{1,n}$, if
\[
(x,y)\in\rho\longrightarrow (u[\bar{w}\,|_ix],
u[\bar{w}\,|_iy])\in\rho
\]
for all $x,y,u\in G$, $\bar{w}\in G^n$,
\item \textit{$v$-negative}, if
\[
 (x[y_1\ldots y_n],y_i)\in\rho
\]
for all $x,y_1,\ldots,y_n\in G$ and $i\in\overline{1,n}$.
\end{itemize}
On $\mathcal{G}$ we define two binary relations $\zeta$ and $\chi$
putting
\[
(x,y)\in\zeta\longleftrightarrow x=y[R_1x\ldots R_nx], \qquad
(x,y)\in\chi \longleftrightarrow (R_1x,R_1y)\in\zeta .
\]
The first relation is a stable order, the second is an $l$-regular
and $v$-negative quasi-order containing $\zeta$ (see \cite{Tr-1}).
For these two relations the following conditions are valid:
\[
\begin{array}{lll}
x\leqslant y\longrightarrow R_ix\leqslant R_iy, \ \ \
i\in\overline{1,n}, &&x\sqsubset y\longleftrightarrow
R_ix\leqslant R_iy, \ \ \ i\in\overline{1,n},
\\[4pt]
x\sqsubset y\longleftrightarrow x[R_1y\ldots R_ny]=x,&&
(R_ix)[y_1\ldots y_n]\leqslant y_i, \ \ \ i\in\overline{1,n},
\\[4pt]
x[R_1y_1\ldots R_ny_n]\leqslant x, &&R_ix=R_iR_kx, \ \ \
i,k\in\overline{1,n},
\end{array}
\]
where $x\leqslant y\longleftrightarrow (x,y)\in\zeta$, and
$x\sqsubset y\longleftrightarrow (x,y)\in\chi$.

Let $W$ be the empty set or an $l$-ideal which is an
$\mathcal{E}$-class of a $v$-regular equivalence relation
$\mathcal{E}$ defined on a Menger algebra $(G,o)$ of rank $n$.
Denote by $(H_a)_{a\in A_{\mathcal{E}}}$ the family of all
$\mathcal{E}$-classes (uniquely indexed by elements of some set
$A_{\mathcal{E}}$) such that $H_a\ne W$. Next, for every $g\in G$
we define on $A_{\mathcal{E}}$ an $n$-place function
$P_{(\mathcal{E},W)}(g)$ putting
\begin{equation} \label{simplrep}
P_{(\mathcal{E},W)}(g)(a_1,\ldots,a_n)=b\longleftrightarrow
g[H_{a_1}\ldots H_{a_n}]\subset H_b,
\end{equation}
where $(a_1,\ldots,a_n)\in\pr
P_{(\mathcal{E},W)}(g)\longleftrightarrow g[H_{a_1}\ldots
H_{a_n}]\cap W=\emptyset$, and $H_b$ is an $\mathcal{E}$-class
containing all elements of the form $g[h_1\ldots h_n]$,\, $h_i\in
H_{a_i}$, $i\in\overline{1,n}$. It is not difficult to see that
the map $P_{(\mathcal{E},W)}\colon\,g\mapsto
P_{(\mathcal{E},W)}(g)$ satisfies the identity
\begin{equation}\label{srep}
  P_{(\mathcal{E},W)}(g[g_1\ldots g_n])=
P_{(\mathcal{E},W)}(g)[P_{(\mathcal{E},W)}(g_1)\ldots
P_{(\mathcal{E},W)}(g_n)],
\end{equation}
i.e., $P_{(\mathcal{E},W)}$ is a representation of $(G,o)$ by
$n$-place functions. This representation will be called
\textit{simplest}.

\section{Stabilizers}

\begin{Definition}\label{D-1}\rm
A nonempty subset $H$ of $G$ is called a \textit{stabilizer} of a
functional Menger system $\mathcal{G}$ (or a functional Menger
$\curlywedge$-algebra
$\mathcal{G}_{\curlywedge}=(G,\curlywedge,o,R_1,\ldots,R_n)$) of
rank $n$ if there exists a representation $P$ of $\mathcal{G}$
(respectively, $\mathcal{G}_{\curlywedge}$) by $n$-place functions
on some set $A$, such that $H=H^a_P$ for some point $a\in A$
common for all elements from $H$.
\end{Definition}

\begin{Theorem}\label{T-1}
For a nonempty subset $H$ of $G$ to be a stabilizer of a
functional Menger system $\,\mathcal{G}$ of rank $n$, it is
necessary and sufficient to be a quasi-stable $l$-unitary normal
$v$-complex contained in some subset $U$ of $G$ such that
$R_iU\subset H,$ \ $R_i(G\!\setminus\! U)\subset G\!\setminus\!U$
and
\begin{eqnarray}
 \label{f-1}
x\in H\wedge y\in H\wedge\,t(x)\in U\longrightarrow t(y)\in U, \\[4pt]
 \label{f-2}
x=y[R_1x\ldots R_nx]\in U\wedge\,u[\bar{w}\,|_iy]\in H
\longrightarrow u[\bar{w}\,|_ix]\in H, \\[4pt]
 \label{f-3}
x=y[R_1x\ldots R_nx]\in U\wedge\,u[\bar{w}\,|_iy]\in U
\longrightarrow u[\bar{w}\,|_ix]\in U
\end{eqnarray}
for all $x,y\in G$, $\bar{w}\in G^{\,n},$ $t\in T_n(G)$ and
$i\in\overline{1,n}$, where the symbol $u[\bar{w}\,|_i\ ]$ may be
empty.\footnote{If $u[\bar{w}\,|_i\ ]$ is the empty symbol, then
$u[\bar{w}\,|_ix]$ is equal to $x$.}
\end{Theorem}
\begin{proof} \textit{Necessity}. Let $H^a_{\Phi}$ be a
stabilizer of $a$ for a functional Menger system
$(\Phi,\mathcal{O},\mathcal{R}_1,\ldots,\mathcal{R}_n)$ of
$n$-place functions. If $f\in H^a_{\Phi}$, i.e.,
$f(a,\ldots,a)=a$, then
\[
f[f\ldots f](a,\ldots,a)=f(f(a,\ldots,a),\ldots,f(a,\ldots,a))=
f(a,\ldots,a)=a.
\]
Thus $f[f\ldots f]\in H^a_{\Phi}$. This means that $H^a_{\Phi}$ is
quasi-stable.

Since, for $f[g\ldots g]\in H^a_{\Phi}$ and $g\in H^a_{\Phi}$ we
have
\[
a=f[g\ldots g](a,\ldots,a)=f(g(a,\ldots,a),\ldots,g(a,\ldots,a))=
f(a,\ldots,a),
\]
then $f\in H^a_{\Phi}$, therefore $H^a_{\Phi}$ is $l$-unitary.

Moreover, if $f(\bar{a})=g(\bar{a})$ for some $f,g\in\Phi$, where
$\bar{a}=(a_1,\ldots,a_n)$, then, as it is not difficult to see,
$t(f)(\bar{a})=t(g)(\bar{a})$ for every $t\in T_n(\Phi)$. So,
$f,g,t (f)\in H^a_{\Phi}$ implies $t(g)\in H^a_{\Phi}$. Therefore
$H^a_{\Phi}$ is a normal $v$-complex.

It is clear that $H^a_{\Phi}\subset
U^a_{\Phi}=\{f\in\Phi\,|\,(a,\ldots,a)\in\pr f\}$,
$\mathcal{R}_iU^a_{\Phi}\subset H^a_{\Phi}$ and $\mathcal{R}_i
(\Phi\!\setminus\!U^a_{\Phi})\subset\Phi\!\setminus\!U^a_{\Phi}$.
If $t(f)\in U^a_{\Phi}$ for some $f,g\in H^a_{\Phi}$ and $t\in
T_n(\Phi)$, then $(a,\ldots,a)\in\pr t(f)$, which, together with
$f(a,\ldots,a)=a=g(a,\ldots,a)$, gives $t(f)(a,\ldots,a)=
t(g)(a,\ldots,a)$. Therefore, $t(g)\in U^a_{\Phi}$. So, the
condition (\ref{f-1}) is satisfied.

To prove (\ref{f-2}) assume
$f=g[\mathcal{R}_1f\ldots\mathcal{R}_nf]\in U^a_{\Phi}$, i.e.,
$f\subset g$ and $f\in U^a_{\Phi}$. This implies
$(a,\ldots,a)\in\pr f$ and $f(a,\ldots,a)=g(a,\ldots,a)$. So, for
$\alpha [\bar{\chi}\,|_ig]\in H^a_{\Phi}$, where $\alpha\in\Phi$
and $\bar{\chi}=(\chi_1,\ldots,\chi_n)\in\Phi^n$, we have
\[\arraycolsep=.5mm\begin{array}{rcl}
a&=&
\alpha[\bar{\chi}\,|_ig](a,\ldots,a)=\alpha(\bar{\chi}(a,\ldots,a)
\,|_ig (a,\ldots,a))\\[4pt]
&=& \alpha(\bar{\chi}(a,\ldots,a)\,|_if(a,\ldots,a))=
\alpha[\bar{\chi}\,|_if](a,\ldots,a),
\end{array}
\]
where $\bar{\chi}(a,\ldots,a)$ is
$\chi_1(a,\ldots,a),\ldots,\chi_n(a,\ldots,a)$. Thus
$\alpha[\bar{\chi}\,|_if]\in H^a_{\Phi}$, which completes the
proof of (\ref{f-2}). The proof of (\ref{f-3}) is analogous.

\medskip
\textit{Sufficiency}. Let $H$ and $U$ be two subsets of $G$
satisfying all the conditions of the theorem. First we shall prove
the following implications:
\begin{eqnarray}
\label{f-4} x\leqslant y\wedge x\in H\longrightarrow y\in H, \\[4pt]
\label{f-5} x\sqsubset y\wedge x\in U\longrightarrow y\in U.
\end{eqnarray}
Indeed, $x\leqslant y$ means $x=y[R_1x\ldots R_nx]$. Since $x\in
H,$ $H\subset U$ and $R_iU\subset H$, we have $R_ix\in H$ for
every $i\in\overline{1,n}$. This, together with the fact that $H$
is an $l$-unitary normal $v$-complex, implies that $H$ is
$v$-unitary. So, $y[R_1x\ldots R_nx]\in H$ and $R_ix\in H$ for
every $i\in\overline{1,n}$, which by the $v$-unitarity of $H$
gives $y\in H$. This proves (\ref{f-4}).

Now, if $x\sqsubset y$ and $x\in U$, then $R_1x\leqslant R_1y$,
i.e., $R_1x=(R_1y)[R_1x\ldots R_nx]$. From $R_iU\subset H$ it
follows $R_ix\in H$, so, applying the $v$-unitarity of $H$ to
$(R_1y)[R_1x\ldots R_nx]\in H$ we obtain $R_1y\in H$, whence we
get $R_1y\in U$. Since $R_i(G\!\setminus\!U)\subset
G\!\setminus\!U$ means that
\begin{equation} \label{f-6}
  R_ix\in U \longrightarrow x\in U
\end{equation}
for every $x\in G$ and $i\in\overline{1,n}$, from $R_1y\in U $ it
follows $y\in U$. This completes the proof of (\ref{f-5}).

The set $G\!\setminus\!U$ is an $l$-ideal of $\mathcal{G}$. In
fact, the $v$-negativity of $\chi$ and $x\in G\!\setminus\!U$
imply $u[\bar{w}\,|_ix]\sqsubset x$, whence, according to
(\ref{f-5}), we conclude $u[\bar{w}\,|_ix]\in G\!\setminus\!U$ for
all $i\in\overline{1,n}$, $u\in G$ and $\bar{w}\in G^{\,n}$. So,
$G\!\setminus\!U$ is an $l$-ideal. Using this fact it is easy to
show that
\begin{eqnarray}
\label{f-7} x\leqslant y\wedge x\in U\wedge\,t(y)\in H
\longrightarrow t(x)\in H, \\[4pt]
\label{f-8} x\leqslant y\wedge x\in U\wedge\,t(y)\in U
\longrightarrow t(x)\in U
\end{eqnarray}
for all $x,y\in G$ and $t\in T_n(G)$.

On $G$ we define two binary relations $ \mathcal{E}_H$ and
$\mathcal{E}_U$ putting
\begin{eqnarray}
 \mathcal{E}_H= \left\{(x,y)\,|\,(\forall t\in T_n(G))
\Big(t(x)\in H\longleftrightarrow t(y)\in H\Big)\right\},
\nonumber \\[4pt]
\mathcal{E}_U=\left\{(x,y)\,|\,(\forall t\in T_n(G))\Big(t(x)\in U
\longleftrightarrow t(y)\in U\Big)\right\}. \nonumber
\end{eqnarray}
These relations are $v$-regular equivalences.
$\mathcal{E}=\mathcal{E}_H\cap\mathcal{E}_U$ also is a $v$-regular
equivalence. For this relation we have
\begin{equation} \label{f-9}
x[\mathcal{E}\langle y_1\rangle\ldots\mathcal{E}\langle y_n
\rangle]\subset\mathcal{E}\langle x[y_1\ldots y_n]\rangle
\end{equation}
for all $x,y_1,\ldots,y_n\in G$, where $\mathcal{E}\langle
y_i\rangle$ denotes an equivalence class of $\mathcal{E}$
containing $y_i$. Moreover, $G\!\setminus\!U$ is an equivalence
class of $\mathcal{E}$.

Also $H$ is an $\mathcal{E}$-class. To prove this fact it is
sufficient to verify the following two conditions:
\begin{eqnarray}
\label{f-10} g_1\in H\wedge\,g_2\in H\longrightarrow g_1
\equiv g_2(\mathcal{E}), \\[4pt]
\label{f-11} g_1\equiv g_2(\mathcal{E})\wedge\,g_1\in H
\longrightarrow g_2\in H.
\end{eqnarray}

Let $g_1,g_2\in H$ and $t(g_1)\in U$ for some $t\in T_n(G)$. Then,
from (\ref{f-1}), it follows $t(g_2)\in U$. Similarly, from
$g_1,g_2\in H$ and $t(g_2)\in U$, we conclude $t(g_1)\in U$.
Hence, $g_1\equiv g_2(\mathcal{E}_U)$. If $g_1,g_2,t(g_1)\in H$,
then, in view of the fact that $H$ is a normal $v$-complex, we
have $t(g_2)\in H$. Similarly, from $g_1,g_2,t(g_2)\in H$ we
deduce $t(g_1)\in H$. So, $g_1\equiv g_2(\mathcal{E}_H)$. Thus
$g_1\equiv g_2(\mathcal{E})$. This proves (\ref{f-10}).

Now let $g_1\in H$ and $g_1\equiv g_2(\mathcal{E})$. Then
$g_1\equiv g_2(\mathcal{E}_H)$, which means that for all $t\in
T_n(G)$ we have $t(g_1)\in H\longleftrightarrow t(g_2)\in H$,
whence $g_1\in H\longleftrightarrow g_2\in H$. So, $g_2\in H$,
i.e., (\ref{f-11}) is proved. Consequently, $H$ is an
$\mathcal{E}$-class.

Let $W=G\!\setminus\!U$. For every $g\in G$ we consider an
$n$-place function $P_{(\mathcal{E},W)}(g)$ defined by
(\ref{simplrep}). Since the map
$P_{(\mathcal{E},W)}\colon\,g\mapsto P_{(\mathcal{E},W)}(g)$
satisfies (\ref{srep}), it is a homomorphism with respect to the
operation $o$. It satisfies also the identity
\begin{equation} \label{f-15}
 P_{(\mathcal{E},W)}(R_ig)=\mathcal{R}_iP_{(\mathcal{E},W)}(g).
\end{equation}
Indeed, for every
$$
\bar{a}=(a_1,\ldots,a_n)\in\pr P_{(\mathcal{E},W)}(R_ig),
$$
where $H_{a_i}=\mathcal{E}\langle x_i\rangle$, $x_i\in U$,
$i\in\overline{1,n}$, we have $(R_ig)[x_1\ldots x_n]\in U$.
Whence, according to $\mathbf{A}_7$, we obtain $x_i[R_1g[x_1\ldots
x_n]\ldots R_ng[x_1\ldots x_n]]\in U$. As $G\!\setminus\!U$ is an
$l$-ideal, the last condition implies $R_ig[x_1\ldots x_n]\in U$
for $i\in\overline{1,n}$. Thus $g[x_1\ldots x_n]\in U$, because
$g[x_1\ldots x_n]\not\in U$ implies $x_i\in G\!\setminus\!U$. So,
$\bar{a}\in\pr P_{(\mathcal{E},W)}(g)$. This proves $\pr
P_{(\mathcal{E},W)}(R_ig)\subset\pr P_{(\mathcal{E},W)}(g)$.

To prove the converse inclusion let $\bar{a}\in\pr
P_{(\mathcal{E},W)}(g),$ where $H_{a_i}=\mathcal{E}\langle
x_i\rangle$, $x_i\in U$ for $i\in\overline{1,n}$. Then
$g[x_1\ldots x_n]\in U$, which, by $R_kU\subset H\subset U$, gives
$R_kg[x_1\ldots x_n]\in U$. Whence, by $\mathbf{A}_6$, we get
$R_k(R_i)[x_1\ldots x_n]\in U$ for $k,i\in\overline{1,n}$. From
this, in view of (\ref{f-6}), we deduce $(R_ig)[x_1\ldots x_n]\in
U$. Thus $\bar{a}\in\pr P_{(\mathcal{E},W)}(R_ig)$ for every
$i\in\overline{1,n}.$

In this way we have proved
\begin{equation} \label{f-16}
\pr P_{(\mathcal{E},W)}(R_ig)=\pr P_{(\mathcal{E},W)}(g)= \pr
\mathcal{R}_iP_{(\mathcal{E},W)}(g)
\end{equation}
for every $g\in G$.

Let $\bar{a}\in\pr P_{(\mathcal{E},W)}(R_ig)$, i.e.,
$(R_ig)[x_1\ldots x_n]\in U$, where $H_{a_i}=\mathcal{E}\langle
x_i\rangle$, $x_i\in U$, $i\in\overline{1,n}$. Applying the
stability of $\zeta$ to $(R_ig)[x_1\ldots x_n]\leqslant x_i$ we
obtain $t((R_ig)[x_1\ldots x_n])\leqslant t(x_i)$ for every $t\in
T_n(G)$ and $i\in\overline{1,n}$. If $t((R_ig)[x_1\ldots x_n])\in
H$, then, according (\ref{f-4}), we get $t(x_i)\in H$. For
$t(x_i)\in H$, in view of (\ref{f-7}), from $(R_ig)[x_1\ldots
x_n]\leqslant x_i$ and $(R_ig) [x_1\ldots x_n]\in U$ we deduce
$t((R_ig)[x_1\ldots x_n])\in H$. So, $(R_ig)[x_1\ldots x_n]\equiv
x_i(\mathcal{E}_H)$. Similarly we can prove $(R_ig)[x_1\ldots
x_n]\equiv x_i(\mathcal{E}_U)$. Thus $(R_ig) [x_1\ldots x_n]
\equiv x_i (\mathcal{E})$ for every $i\in\overline{1,n}$.
Therefore $P_{(\mathcal{E},W)}(R_ig)(a_1,\ldots,a_n)=a_i$.
Consequently,
 \[
 P_{(\mathcal{E},W)}(R_ig)(a_1,\ldots,a_n)=
(\mathcal{R}_iP_{(\mathcal{E},W)}(g))(a_1,\ldots,a_n).
\]
This, together with (\ref{f-16}), gives
(\ref{f-15}).

From (\ref{srep}) and (\ref{f-15}) it follows that
$P_{(\mathcal{E},W)}$ is a representation of $\mathcal{G}$ by
$n$-place functions.

Observe that
\begin{equation} \label{f-18}
  g\in H\longleftrightarrow P_{(\mathcal{E},W)}(g)(a,\ldots,a)=a,
\end{equation}
where $a$ is this element of $A_{\mathcal{E}}$ which is used as
index of the $\mathcal{E}$-class $H$. In fact, for $g\in H$ the
quasi-stability of $H$ implies $g[g\ldots g]\in H$. Whence,
$g[H\ldots H]\subset H$ because $H$ is an $\mathcal{E}$-class and
the relation $\mathcal{E}$ is $v$-regular. So,
$P_{(\mathcal{E},W)}(g)(a,\ldots,a)=a$. Conversely, if $g[H\ldots
H]\subset H$, then $g[h\ldots h]\in H$ for every $h\in H$. From
this, by the $l$-unitarity of $H$, we get $g\in H$, which
completes the proof of (\ref{f-18}).

To complete this proof we remind that any algebra $\mathcal{G}$
satisfying the axioms $\mathbf{A}_1-\mathbf{A}_7$, has a faithful
representation by $n$-place functions \cite{Tr-1}. Let $P_1$ be
this representation. Then, as it is not difficult to verify,
$P=P_1+P_{(\mathcal{E},W)}$ is a faithful representation of
$\mathcal{G}$ for which $H=H^a_P$. This completes the proof.
\end{proof}

Let $\mathcal{G}$ be a functional Menger system of rank $n$ and
$H$ be some subset of $G$. We say that a subset $X$ of $G $ is
\textit{$C_H$-closed}, if for all $a,b,c\in G,$ $t\in T_n(G)$ the
implication:
 \[
\left.
  \begin{array}{l}
    a=b [R_1a\ldots R_na]\vee a,b\in H, \\[4pt]
    t(a)[R_1c\ldots R_nc]=t(a), \\[4pt]
    a,t(b)\in X \
  \end{array}\right \}\longrightarrow c\in X
\]
is valid. In the abbreviated form this implication can be written
as
\begin{equation} \label{f-20}
  (a\leqslant b\vee a,b\in H)\,\wedge\, t(a)\sqsubset c\,\wedge\,
  a,t(b)\in X\longrightarrow c\in X.
\end{equation}

Let $C_H(X)$ denotes the set of all $c\in G$ for which there exist
$a,b\in G$ and $t\in T_n(G)$ such that the premise of (\ref{f-20})
is satisfied. Further, let
\[
 C_H[X]=\bigcup\limits_{m=0}^{\infty}\stackrel{m}{C}_H\!(X),
\]
where $\stackrel{0}{C}_H\!(X)=X$,
$\,\stackrel{m+1}{C_H}\!\!(X)=C_H(\stackrel{m}{C}_H\!(X))$ for
every $m=0,1,2,\ldots$

By induction we can prove that $g\in\,\stackrel{m}{C}_H\!\!(X)$ if
and only if the following system of conditions is fulfilled:
\begin{equation}\label{f-21}
 \left.
\begin{array}{l}
(\,a_1=b_1 [R_1a\ldots R_na]\vee a_1,b_1\in H\,)\,\wedge\, t_1(a_1)[R_1g\ldots R_ng]=t_1(a_1) \\[4pt]
 \bigwedge\limits_{i=1}^{2^{m-1}-1} \left(\begin{array}{l}
    a_{2i}=b_{2i}[R_1a_{2i}\ldots R_na_{2i}]\vee\, a_{2i},b_{2i}\in H, \\[4pt]
    t_{2i}(a_{2i})[R_1a_i\ldots R_na_i]=t_{2i}(a_{2i}), \\[4pt]
    a_{2i+1}=b_{2i+1}[R_1a_{2i+1}\ldots R_na_{2i+1}]\vee a_{2i+1},b_{2i+1}\in H, \\[4pt]
    t_{2i+1}(a_{2i+1})[R_1t_i(b_i)\ldots R_nt_i(b_i)]=t_{2i+1}(a_{2i+1}) \
  \end{array}\right) \\[4pt]
  \bigwedge\limits_{i=2^{m-1}}^{2^m-1}(a_i\in X\wedge\,t_i(b_i)\in X)
\end{array}\right\}
\end{equation}
where $a_k,b_k\in G$, $t_k\in T_n(G)$.

In the sequel, the system of conditions (\ref{f-21}) will be
denoted by $\mathfrak{M}_H(X,m,g)$.

\begin{Theorem}\label{T-2}
Let $\mathcal{G}$ be a functional Menger system of rank $n$. A
nonempty subset $H$ of $G$ is a stabilizer of $\mathcal{G}$ if and
only if it is a quasi-stable $l$-unitary normal $v$-complex such
that $R_iH\subset H$ for every $i\in\overline{1,n}$ and
\begin{equation} \label{f-22}
 x=y[R_1x\ldots R_nx]\in C_H[H]\,\wedge\,u[\bar{w}\,|_iy]\in H
 \longrightarrow u[\bar{w}\,|_ix]\in H
\end{equation}
for all $x,y,u\in G,$ $\bar{w}\in G^{\,n}$, $i\in\overline{1,n}$,
where the symbol $u[\bar{w}\,|_i\ ]$ may be empty.
\end{Theorem}
\begin{proof} \textit{Necessity}. Let $H^a_{\Phi}$ be the
stabilizer of a point $a$ in a functional Menger system
$(\Phi,\mathcal{O},\mathcal{R}_1,\ldots,\mathcal{R}_n)$ of
$n$-place functions. Obviously it is a quasi-stable $l$-unitary
normal $v$-complex of
$(\Phi,\mathcal{O},\mathcal{R}_1,\ldots,\mathcal{R}_n)$. If $f\in
H^a_{\Phi}$, then $f(a,\ldots,a)=a$, whence
$\mathcal{R}_if(a,\ldots,a)=a$, i.e., $\mathcal{R}_if\in
H^a_{\Phi}$. So, $\mathcal{R}_iH^a_{\Phi}\subset H^a_{\Phi}$ for
every $i\in\overline{1,n}$.

To prove (\ref{f-22}), we shall consider
$f=g[\mathcal{R}_1f\ldots\mathcal{R}_nf]\in
C_{H^a_{\Phi}}[H^a_{\Phi}]$ and $\alpha [\bar{\omega}\,|_ig]\in
H^a_{\Phi}$ for some $f,g,\alpha\in\Phi$,
$\bar{\omega}\in\Phi^{n}$, $i\in\overline{1,n}$, where
$C_{H^a_{\Phi}}[H^a_{\Phi}]=
\bigcup\limits_{m=0}^{\infty}\stackrel{m}{C}_{H^a_{\Phi}}\!(H^a_{\Phi})$.
Then
$f=g[\mathcal{R}_1f\ldots\mathcal{R}_nf]\in\,\stackrel{m}{C}_{H^a_{\Phi}}\!(H^a_{\Phi})$
for some $m\in\mathbb {N}$. But, as it is easily to see by
induction,
$\varphi\in\,\stackrel{m}{C}_{H^a_{\Phi}}\!(H^a_{\Phi})$ implies
$(a,\ldots,a)\in\pr\varphi$. The above means that
$(a,\ldots,a)\in\pr f$ and $(a,\ldots,a)\in\pr
g[\mathcal{R}_1f\ldots\mathcal{R}_nf]$. Therefore
 \[\arraycolsep=.5mm
\begin{array}{rcl}
 f(a,\ldots,a)&=&g[\mathcal{R}_1f\ldots\mathcal{R}_nf](a,\ldots,a)\\[4pt]
&=&g(\mathcal{R}_1f(a,\ldots,a),\ldots,\mathcal{R}_nf(a,\ldots,a))=
g(a,\ldots,a).
\end{array}
\]
Moreover, for $\alpha[\bar{\omega} \,|_ig]\in H^a_{\Phi}$ we have
\[\arraycolsep=.5mm
\begin{array}{rcl}
a&=&\alpha [\bar{\omega}\,|_ig](a,\ldots,a)=
\alpha(\bar{\omega}(a,\ldots,a)\,|_ig (a,\ldots,a))\\[4pt]
&=&\alpha(\bar{\omega}(a,\ldots,a)\,|_if(a,\ldots,a))=
\alpha[\bar{\omega}\,|_if](a,\ldots,a),
\end{array}\]
where $\bar{\omega}(a,\ldots,a)$ denotes
$\omega_1(a,\ldots,a),\ldots,\omega_n(a,\ldots,a)$. So,
$\alpha[\bar{\omega}\,|_if]\in H^a_{\Phi}$, which completes the
proof of (\ref{f-22}).

\textit{Sufficiency}. Let $H$ satisfy all the conditions of the
theorem. We prove that it satisfies also all the conditions of
Theorem~\ref{T-1}. Since by the assumption $H$ is a quasi-stable
$l$-unitary normal $v$-complex such that $R_iH\subset H$ for every
$i\in\overline{1,n}$, we must prove that it satisfies the
conditions (\ref{f-1}), (\ref{f-2}), (\ref{f-3}) and $R_iU\subset
H$, $R_i(G\!\setminus\!U)\subset G\!\setminus\!U$ for some
$U\subset G$ containing $H$ and all $i\in\overline{1,n}$.

First we prove that $H$ satisfies the condition (\ref{f-4}).
Indeed, if $x\leqslant y $ and $x\in H$, then $y[R_1x\ldots
R_nx]=x\in H$ and $R_ix\in H,$ $\,i\in\overline{1,n}$. As $H$ is
an $l$-unitary normal $v$-complex, $H$ is a $v$-unitary subset.
Therefore $y\in H,$ which completes the proof of (\ref{f-4}).

Now let $U=C_H[H]$. Clearly $H\subset U$. Since $C_H[H]=
\bigcup\limits_{m=0}^{\infty}\!\stackrel{m}{C}_H\!(H)$, to prove
that $R_iU\subset H$ for every $i\in\overline{1,n},$ it is
sufficient to show that for every $m\in\mathbb {N}$ holds the
inclusion
\begin{equation} \label{f-23}
  R_i(\stackrel{m}{C}_H\!(H))\subset H .
\end{equation}
For $m=0$ it is obvious because $\stackrel{0}{C}_H\!(H)=H$ and
$R_iH\subset H$ for all $i\in\overline{1,n}$. Suppose that
(\ref{f-23}) is valid for some $k\in\mathbb{N}$. We prove that it
is valid for $k+1$. Let $g\in\,\stackrel{k+1}{C_H}\!(H)$. Then
 \[
(a\leqslant b\vee a,b\in H)\wedge\,t(a)\sqsubset g \wedge\,
a,t(b)\in\,\stackrel{k}{C}_H\!(H)
\]
for some $a,b\in G$, $t\in T_n(G)$. From
$a,t(b)\in\,\stackrel{k}{C}\!(H)$, according to our supposition,
we get $R_ia,R_it(b)\in H$. If $a\leqslant b$ and $a\in H$, then,
by (\ref{f-4}), we have $b\in H$. Thus $a,b,R_it(b)\in H$. Whence
$R_it(a)\in H$ because $H$ is a normal $v$-complex. The condition
$t(a)\sqsubset g$ implies $R_it(a)\leqslant R_ig$, which, by
(\ref{f-4}), gives $R_ig\in H$. In a similar way $a,b\in H$ and
$t(a)\sqsubset g$ proves $R_ig\in H$. Thus we have shown that
$R_i(\stackrel{k+1}{C_H}\!(H))\subset H$. So, the inclusion
(\ref{f-23}) is valid for $m=k+1$ and consequently for every
$m\in\mathbb{N}$. Therefore $R_iU\subset H$ for every
$i\in\overline{1,n}.$

To prove the inclusion $R_i(G\!\setminus\!U)\subset
G\!\setminus\!U$ observe that it is equivalent to the condition
\[
(\forall g\in G)\,(g\in G\!\setminus\!U\longrightarrow R_ig\in
G\!\setminus\!U),
\]
which can be written in the form
\begin{equation} \label{f-24}
  (\forall g\in G)\,(R_ig\in U\longrightarrow g\in U).
\end{equation}
In the case $U=C_H[H]$ the last condition means that
\begin{equation} \label{f-25}
(\forall g\in G)(\forall m\in\mathbb{N})\,
\left(R_ig\in\,\stackrel{m}{C}_H \!(H)\longrightarrow (\exists
n\in\mathbb{N})\, g\in\,\stackrel{n}{C}_H\!(H)\right).
\end{equation}
Let $R_ig\in\,\stackrel{m}{C}_H\!(H)$ for some $g\in G$ and
$m\in\mathbb{N}$. Considering $R_ig=R_i(R_ig)$ and (\ref{f-23}),
we conclude $R_ig\in H$. Thus $R_ig\leqslant R_ig$, $R_ig\sqsubset
g$ and $R_ig\in H$. Therefore $g\in C_H(H)$, i.e.
$g\in\,\stackrel{n}{C}_H\!(H)$ for some $n\in\mathbb{N}$. This
proves (\ref{f-25}). So, $R_i(G\!\setminus\!U)\subset
G\!\setminus\!U$ for $U=C_H[H]$ and $i\in\overline{1,n}$.

To prove (\ref{f-1}) observe that for $x,y\in H$ and $t(x)\in
U=C_H[H]$, just proved inclusion implies $R_it(x)\in R_iU\subset
H$. Thus $x,y,R_it(x)\in H$. But $H$ is a normal $v$-complex,
hence $R_it(y)\in H$. So, for $U=C[H]$, the condition (\ref{f-1})
is satisfied. Also (\ref{f-2}) is valid because for $U=C_H[H]$ it
coincides with (\ref{f-22}). Since the subset $C_H[H]$ is
$C_H$-closed, the condition (\ref{f-3}) is valid too. This means
that $H$ satisfies all the conditions of Theorem~\ref{T-1}. Hence,
$H$ is a stabilizer of a functional Menger system $\mathcal{G}$.
 \end{proof}

As it is not difficult to see, the condition (\ref{f-22}) is
equivalent to the system of conditions $(A'_m)_{m\in\mathbb {N}}$,
where
 \[
A'_m\,\colon\ x=y [R_1x\ldots R_nx]\,\wedge\,
x\in\,\stackrel{m}{C}_H\!(H)\,\wedge\, u[\bar{w}\,|_iy]\in H
\longrightarrow u[\bar{w}\,|_ix]\in H
 \]
for all $x,y,u\in G$, $\bar{w}\in G^n$, $i\in\overline{1,n}$.
Since, by (\ref{f-21}), the condition $A'_m $ is equivalent to
 \[
A_m\,\colon\ x=y[R_1x\ldots
R_nx]\,\wedge\,\mathfrak{M}_H(H,m,x)\,\wedge\, u[\bar{w}\,|_iy]\in
H \longrightarrow u[\bar{w}\,|_ix]\in H,
 \]
the last theorem can be written in the form:

\begin{Theorem} \label{T-3}
A nonempty subset $H$ of $G$ is a stabilizer of a functional
Menger system $\,\mathcal{G}$ of rank $n$ if and only if it is a
quasi-stable $l$-unitary normal $v$-complex such that $R_iH\subset
H$ for every $i\in\overline{1,n}$ and the system of conditions
$(A_m)_{m\in\mathbb {N}}$ is satisfied.
\end{Theorem}

Now we shall characterize stabilizers of  functional Menger
$\curlywedge$-algebras.

\begin{Theorem} \label{T-4}
A nonempty subset $H$ of $G$ is a stabilizer of a functional
Menger $\curlywedge$-algebra
$\mathcal{G}_{\curlywedge}=(G,o,\curlywedge,R_1,\ldots,R_n)$ of
rank $n$ if and only if
\begin{enumerate}
\item[$1)$] it is a quasi-stable, $\curlywedge$-stable and $v$-unitary
subset of $G$,
\item[$2)$] there exists a subset $U$ of $G$ such that $H\subset U$,
$R_iU\subset H$ and $R_i(G\!\setminus\!U)\subset G\!\setminus\!U$
for every $i\in\overline{1,n}$,
\item[$3)$] the following two implications
\begin{eqnarray}
 \label{T4-1} x\in U\wedge y\in H\longrightarrow y[R_1x\ldots R_nx]\in H, \\[4pt]
 \label{T4-1a} x\in U \wedge y\in U\longrightarrow y[R_1x\ldots R_nx]\in U
\end{eqnarray}
are valid for all $x,y\in G$.
\end{enumerate}
\end{Theorem}
\begin{proof} The proof of the necessity of the conditions of the theorem is
similar to the proof of the necessity of the conditions of
Theorem~\ref{T-1}. So, we prove only their sufficiency.

Let all these conditions be satisfied. First we shall show the
implication
\begin{equation}
 \label{T4-2} x\curlywedge y\in U\wedge\, u[\bar{w}\,|_i(y\curlywedge
z)]\in U\longrightarrow u[\bar{w}\,|_i(x\curlywedge y\curlywedge
z)] \in U
\end{equation}
for $x,y,z\in G$ and $i\in\overline{1,n}$, where the symbol
$u[\bar{w}\,|_i\ ]$ may be empty. For this suppose that the
premise of (\ref{T4-2}) is satisfied. Then, according to
(\ref{T4-1a}), we have
\[
u[\bar{w}\,|_i(y\curlywedge z)][R_1(x\curlywedge y)\ldots
R_n(x\curlywedge y)]\in U .
\]
Whence, by $\mathbf{A}_3$, we obtain $u[\bar{w}\,|_i(y\curlywedge
z)[R_1(x\curlywedge y)\ldots R_n(x\curlywedge y)]]\in U.$ From
this, applying $\mathbf{A}_8$, we conclude
$u[\bar{w}\,|_i(z\curlywedge y[R_1(x\curlywedge y)\ldots
R_n(x\curlywedge y)])]\in U$, which, by $\mathbf{A}_9$, implies
$u[\bar{w}\,|_i(x\curlywedge y\curlywedge z)]\in U$. This
completes the proof of (\ref{T4-2}).

Further, in the same way as in the proof of Theorem~\ref{T-1}, we
can show that the conditions (\ref{f-4}), (\ref{f-5}), (\ref{f-6})
are satisfied and $G\!\setminus\!U$ is an $l$-ideal. Next, we
consider the relation:
\begin{equation} \label{f-28}
\mathcal{E}[U]= \{(g_1,g_2)\in G\times G\,|\,g_1\curlywedge g_2\in
U\,\vee\, g_1,g_2\in U^{\,\prime}\},
\end{equation}
where $U^{\,\prime}=G\!\setminus\!U$. The reflexivity and symmetry
of this relation are obvious. Let
$(g_1,g_2),(g_2,g_3)\in\mathcal{E}[U]$. If $g_1,g_2,g_3\in
U^{\,\prime}$, then evidently $(g_1,g_3)\in\mathcal{E}[U]$. If
$g_1\curlywedge g_2,g_2\curlywedge g_3\in U$, then, according to
(\ref{T4-2}), we have $g_1\curlywedge g_2\curlywedge g_3\in U$.
Whence, in view of (\ref{f-5}) and the fact that $g_1\curlywedge
g_2\curlywedge g_3\leqslant g_1\curlywedge g_3$, \
$\zeta\subset\chi$ and $g_1\curlywedge g_2\curlywedge g_3\sqsubset
g_1\curlywedge g_3$, we conclude $g_1\curlywedge g_3\in U$, i.e.,
$(g_1, g_3)\in\mathcal{E}[U]$. So, $\mathcal{E}[U]$ is also
transitive. Thus $\mathcal{E}[U]$ is an equivalence relation.

It is $v$-regular too. Indeed, if $(g_1,g_2)\in\mathcal{E}[U]$,
then $g_1,g_2\in U^{\,\prime}$ or $g_1\curlywedge g_2\in U$. Since
$U^{\,\prime}$ is an $l$-ideal, in the case $g_1,g_2\in
U^{\,\prime}$ we have $u[\bar{w}\,|_ig_1],u[\bar{w}\,|_ig_2]\in
U^{\,\prime}$, i.e.,
\begin{equation} \label{T4-xx}
 u[\bar{w}\,|_ig_1]\equiv u[\bar{w}\,|_ig_2](\mathcal{E}[U]).
\end{equation}
In the case when $g_1\curlywedge g_2\in U$ and
$u[\bar{w}\,|_ig_1],u[\bar{w}\,|_ig_2]\in U^{\,\prime}$,
(\ref{T4-xx}) is satisfied too. In the case $g_1\curlywedge g_2\in
U$, $u[\bar{w}\,|_ig_1]\in U$, according to (\ref{T4-1a}), we have
$$
u[\bar{w}\,|_ig_1][R_1(g_1\curlywedge g_2)\ldots
R_n(g_1\curlywedge g_2)]\in U,
$$
whence, applying $\mathbf{A}_3$, we obtain
$$
u[\bar{w}\,|_ig_1[R_1(g_1\curlywedge g_2)\ldots R_n(g_1\curlywedge
g_2)]]\in U,
$$
which, by $\mathbf{A}_9$, implies $u[\bar{w}\,|_i(g_1\curlywedge
g_2)]\in U$. This, in view of $u[\bar{w}\,|_i(g_1\curlywedge
g_2)]\sqsubset u[\bar{w}\,|_i g_2]$ and (\ref{f-5}), gives
$u[\bar{w}\,|_i g_2]\in U$.

Similarly we can prove that $g_1\curlywedge g_2\in U$ and
$u[\bar{w}\,|_ig_2]\in U$ imply $u[\bar{w}\,|_i g_1]\in U$.
Therefore $u[\bar{w}\,|_i g_1]$ and $u[\bar{w}\,|_ig_2]$ belong or
do not belong to $U$ simultaneously. Let $u[\bar{w}\,|_i
g_1],u[\bar{w}\,|_i g_2]\in U$. Since
\[
u[\bar{w}\,|_i(g_1\curlywedge g_2)]\leqslant u[\bar{w}\,|_i g_j]
\]
for $j=1,2$, from the above we obtain
\[
u[\bar{w}\,|_i(g_1\curlywedge g_2)]\sqsubset
u[\bar{w}\,|_ig_1]\curlywedge u[\bar{w}\,|_i g_2],
\]
which, after application of (\ref{f-5}), gives
$u[\bar{w}\,|_i(g_1\curlywedge g_2)]\in U$. Thus (\ref{T4-xx}) is
satisfied in any case. So, the relation $\mathcal{E}[U]$ is
$i$-regular for every $i\in\overline{1,n}$, i.e., it is
$v$-regular.

$H$ is an equivalence class of $\mathcal{E}[U]$. Indeed, if
$g_1,g_2\in H$, then $g_1\curlywedge g_2\in H$ by the
$\curlywedge$-stability of $H.$ Consequently $g_1\curlywedge
g_2\in U$, i.e., $g_1\equiv g_2(\mathcal{E}[U])$. On the other
hand, if $g_1\equiv g_2(\mathcal{E}[U])$ and $g_1\in H$, then
$g_1\curlywedge g_2\in U$, whence, by (\ref{T4-1}), we have
\[
g_1[R_1(g_1\curlywedge g_2)\ldots R_n (g_1\curlywedge g_2)]\in H.
\]
From this, applying $\mathbf{A}_9$, we deduce $g_1\curlywedge
g_2\in H$, which, in view of $g_1\curlywedge g_2\leqslant g_2$ and
(\ref{f-4}), implies $g_2\in H$. So, $H$ is an equivalence class
of $\mathcal{E}[U]$.

Consider the simplest representation $P_{(\mathcal{E}[U],
U^{\,\prime})}$ of $(G,o)$ induced by the pair
$(\mathcal{E}[U],U^{\,\prime})$ (see the section 2). We shall
prove that this representation satisfies the following two
identities:
\begin{eqnarray}
\label{pp} P_{(\mathcal{E}[U],U^{\,\prime})}(g_1\curlywedge g_2)=
P_{(\mathcal{E}[U],U^{\,\prime})}(g_1)\cap P_{(\mathcal{E}[U],U^{\,\prime})}(g_2),
\\[4pt]
\label{p} P_{(\mathcal{E}[U],U^{\,\prime})}(R_ig)=
\mathcal{R}_iP_{(\mathcal{E}[U],U^{\,\prime})}(g),
\end{eqnarray}
where $g, g_1,g_2\in G$ and $i\in\overline{1,n}$.

Let $P_{(\mathcal{E}[U],U^{\prime})}(g_1\curlywedge
g_2)(a_1,\ldots, a_n)=b$, where $H_{a_i}$, $H_{b}$ are equivalence
classes of $\mathcal{E}[U]$ containing elements $x_i,y\in U$,
respectively. Thus $(g_1\curlywedge g_2)[x_1\ldots x_n]\equiv
y(\mathcal{E}[U])$, whence
\begin{equation} \label{*}
y\curlywedge g_1 [x_1\ldots x_n]\curlywedge g_2 [x_1\ldots x_n]\in
U ,
\end{equation}
according to $\mathbf{A}_{10}$. Consequently, $y\curlywedge
g_1[x_1\ldots x_n]\in U$ and $y\curlywedge g_2[x_1\ldots x_n]\in
U$, i.e., $g_1[x_1\ldots x_n]\equiv y(\mathcal{E}[U])$ and
$g_2[x_1\ldots x_n]\equiv y(\mathcal{E}[U])$. Therefore
\begin{equation} \label{**}
P_{(\mathcal{E}[U],U^{\,\prime})}(g_1)(a_1,\ldots,a_n)=b \ \
\mbox{ and } \ \
P_{(\mathcal{E}[U],U^{\,\prime})}(g_2)(a_1,\ldots,a_n)=b.
\end{equation}
This proves the inclusion
$P_{(\mathcal{E}[U],U^{\,\prime})}(g_1\curlywedge g_2)\subset
P_{(\mathcal{E}[U],U^{\,\prime})}(g_1)\cap
P_{(\mathcal{E}[U],U^{\,\prime})}(g_2)$.

To prove the inverse inclusion assume (\ref{**}). Then
$g_1[x_1\ldots x_n]\equiv y(\mathcal{E}[U])$ and $g_2[x_1\ldots
x_n]\equiv y(\mathcal{E}[U])$, i.e., $y\curlywedge g_1[x_1\ldots
x_n]\in U$ and $y\curlywedge g_2[x_1\ldots x_n]\in U$. From this,
according to (\ref{T4-2}), we obtain (\ref{*}), whence, as it was
shown above, we have
$P_{(\mathcal{E}[U],U^{\,\prime})}(g_1\curlywedge
g_2)(a_1,\ldots,a_n)=b$. So,
$P_{(\mathcal{E}[U],U^{\,\prime})}(g_1)\cap
P_{(\mathcal{E}[U],U^{\,\prime})}(g_2)\subset
P_{(\mathcal{E}[U],U^{\,\prime})}(g_1\curlywedge g_2)$. This
completes the proof of (\ref{pp}).

Further, using the same method as in the proof of the condition
(\ref{f-16}), we can prove that
\[
\pr P_{(\mathcal{E}[U],U^{\,\prime})}(R_ig)= \pr P_{(\mathcal{E}[U],U^{\,\prime})}(g)=
\pr\mathcal{R}_iP_{(\mathcal{E}[U],U^{\,\prime})}(g)
\]
for every $g\in G.$ Now if $(a_1,\ldots,a_n)\in\pr
P_{(\mathcal{E}[U],U^{\,\prime})}(R_ig)$, then $(R_ig)[x_1\ldots
x_n]\in U$ for all $x_i$ from the class $H_{a_i}$,
$i\in\overline{1,n}$. Since $(R_ig)[x_1\ldots x_n]\leqslant x_i$
implies $(R_ig)[x_1\ldots x_n]\curlywedge x_i= (R_ig)[x_1\ldots
x_n],$ the above gives $(R_ig)[x_1\ldots x_n]\curlywedge x_i\in
U$. Hence $(R_ig)[x_1\ldots x_n]\equiv x_i(\mathcal{E}[U])$.
Therefore,
\[
P_{(\mathcal{E}[U],U^{\,\prime})}(R_ig)(a_1,\ldots,a_n)=a_i
\]
for all $i\in\overline{1,n}$. This proves (\ref{p}). So, the
simplest representation $P_{(\mathcal{E}[U],U')}$ of $(G,o)$ is
also the simplest representation of a functional Menger
$\curlywedge$-algebra $(G,o,\curlywedge,R_1,\ldots,R_n)$.
Moreover, analogously as (\ref{f-18}), we can prove that
\[
g\in H\longleftrightarrow
P_{(\mathcal{E}[U],U')}(g)(a,\ldots,a)=a,
\]
where $a$ is an elements used as index of the
$\mathcal{E}[U]$-class $H.$

Since the algebra $\mathcal{G}_{\curlywedge}$ satisfies the axioms
$\mathbf{A}_1$ - $\mathbf{A}_{10}$, from \cite{Dudtro} it follows
that there exists an isomorphism of $\mathcal{G}$ onto some
functional Menger $\cap$-algebra of $n$-place functions. Denote
this isomorphism by $P_1$ and consider the representation
$P=P_1+P_{(\mathcal{E}[U],U^{\,\prime})}$. It is clear that $P$ is
a faithful representation of $\mathcal{G}$ by $n$-place functions
and $H=H^a_P$.
\end{proof}

\begin{Theorem} \label{T-5}
A nonempty subset $H$ of a functional Menger $\curlywedge$-algebra
$\mathcal{G}_{\curlywedge}$ of rank $n$ is its stabilizer if and
only if it is stable, $\curlywedge$-stable and $v$-unitary subset
of $G$ such that $R_iH\subset H$ for every $i\in\overline{1,n}$.
\end{Theorem}
\begin{proof} The necessity of these conditions is obvious,
therefore we shall prove only their sufficiency.

Let $H$ satisfies all these conditions. Then $\zeta(H)\subset H$
which is equivalent to (\ref{f-4}).\footnote{\,Remind that for any
relation $\rho\subset X\times Y$ and any subset $A$ of $X$ by
$\rho(A)$ is denoted the set $\{y\in Y\,|\,(\exists x\in A)
(x,y)\in\rho\}$.} Indeed, for $x\in H $ and $x\leqslant y$ we have
$y[R_1x\ldots R_nx]=x\in H$ and $R_ix\in H$ for ever
$i\in\overline{1,n}$, whence, according to the $v$-unitarity of
$H$, we obtain $y\in H$. This proves (\ref{f-4}).

Using this condition we shall prove that $U_0=\chi(H)$ and $H$
satisfy all conditions of Theorem~\ref{T-4}. The stability of $H$
implies its quasi-stability. Because $\chi$ is a quasi-order we
have also $H\subset\chi(H)=U_0$. Moreover, for every $x\in U_0$
there exists $h\in H$ such that $(h,x)\in\chi$, i.e.,
$R_ih\leqslant R_ix$ for every $i\in\overline{1,n}$. Since
$R_ih\in H$ for every $h\in H$ and $i\in\overline{1,n}$, the
above, according to (\ref{f-4}), implies $R_ix\in H$ for
$i\in\overline{1,n}$. Therefore $R_iU_0\subset H$ for every
$i\in\overline{1,n}$.

The inclusion $R_i(G\!\setminus\!U_0)\subset G\!\setminus\!U_0$,
where $i\in\overline{1,n}$, is equivalent to the implication
$(\forall x\in G) (x\in G\!\setminus\!U_0\longrightarrow R_ix\in
G\!\setminus\!U_0)$, which, by contraposition, means that
\[
 R_ix\in U_0\longrightarrow x\in U_0
\]
for every $x\in G.$ But $U_0=\chi(H)$, so, for every $x\in G$ and
every $h\in H$
\begin{equation} \label{kr*}
(h,R_ix)\in\chi\longrightarrow x\in\chi (H).
\end{equation}
Let $(h,R_ix)\in\chi$ for some $x\in G$ and $h\in H.$ Then
$$
H\ni h=h[R_1R_ix\ldots R_nR_ix]=h[R_1x\ldots R_nx],
$$
whence $(h,x)\in\chi$. Therefore $x\in\chi(H)$. This means that
the implication (\ref{kr*}) is valid. So,
$R_i(G\!\setminus\!U_0)\subset G\!\setminus\!U_0$ for every
$i\in\overline{1,n}$. In this way we have proved that $U_0$ and
$H$ satisfy the first two conditions of Theorem~\ref{T-4}.

To prove that $U_0$ and $H$ satisfy the third condition of this
theorem, we must show the following two implications:
\begin{equation} \label{v}
x\in U_0 \wedge y\in H \longrightarrow y [R_1x\ldots R_nx]\in H,
\end{equation}
\begin{equation} \label{v*}
x\in U_0\wedge y\in U_0\longrightarrow y[R_1x\ldots R_nx]\in U_0 .
\end{equation}

Let $x\in U_0=\chi(H)$ and $y\in H$. Then, $(h,x)\in\chi$, i.e.,
$h=h[R_1x\ldots R_nx]$ for some $h\in H$. Since from $h\in H$ it
follows $R_ih\in H $ for all $i\in\overline{1,n}$, we have
$y,R_ih\in H$, $i\in\overline{1,n}$. This, by the stability of
$H,$ implies $y[R_1h\ldots R_nh]\in H$. Whence, in view of
\[
\arraycolsep=.5mm\begin{array}{llll}
y[R_1h\ldots R_nh]&=y[R_1h[R_1x\ldots R_nx]\ldots R_nh[R_1x\ldots R_nx]]\\[4pt]
&\stackrel{\mathbf{A}_4}{=}y[(R_1h)[R_1x\ldots R_nx]\ldots (R_nh)[R_1x\ldots R_nx]]\\[4pt]
&\stackrel{\mathbf{A}_1}{=}y[R_1h\ldots R_nh][R_1x\ldots R_nx]\\[4pt]
&\stackrel{\mathbf{A}_5}{=}y[R_1x\ldots R_nx][R_1h\ldots R_nh],
\end{array}
\]
we conclude $y[R_1x\ldots R_nx][R_1h\ldots R_nh]\in H$. But $H$ is
$v$-unitary and $R_ih\in H$ for every $i\in\overline{1,n}$, so,
$y[R_1x\ldots R_nx]\in H$. This proves (\ref{v}).

Now let $x,y\in U_0=\chi(H)$. Then there exist $a,b\in H$ such
that $(a,x)\in\chi$ and $(b,y)\in\chi$, i.e., $a=a[R_1x\ldots
R_nx]$ and $b=b[R_1y\ldots R_ny]$. Because $b\in H$ implies
$R_ib\in H$ for every $i\in\overline{1,n}$, we have $a,R_ib\in H$
for all $i\in\overline{1,n}$. Whence, according to the stability
of $H$, we get $a[R_1b\ldots R_nb]\in H$. Moreover,
\[\arraycolsep=.5mm\begin{array}{llll}
a[R_1b\ldots R_nb]&=a[R_1x\ldots R_nx][R_1b\ldots R_nb] \\[4pt]
&\stackrel{\mathbf{A}_5}{=}a[R_1b\ldots R_nb][R_1x\ldots R_nx] \\[4pt]
&=a[R_1b[R_1y\ldots R_ny]\ldots R_nb [R_1y\ldots R_ny]] \\[4pt]
&\stackrel{\mathbf{A}_4}{=}a[(R_1b)[R_1y\ldots R_ny]\ldots (R_nb)
[R_1y\ldots R_ny]][R_1x\ldots R_nx]\\[4pt]
&\stackrel{\mathbf{A}_1}{=}a[R_1b\ldots R_nb][R_1y\ldots R_ny][R_1x\ldots R_nx]\\[4pt]
&\stackrel{\mathbf{A}_1}{=}a[R_1b\ldots R_nb][(R_1y)[R_1x\ldots
R_nx]\ldots (R_ny)[R_1x\ldots R_nx]]\\[4pt]
&\stackrel{\mathbf{A}_4}{=}a[R_1b\ldots R_nb][R_1(y[R_1x\ldots
R_nx])\ldots R_n(y[R_1x\ldots R_nx])].
\end{array}\]
Thus, $(a [R_1b\ldots R_nb], y [R_1x\ldots R_nx])\in\chi$, i.e.,
$y[R_1x\ldots R_nx]\in U_0=\chi(H)$. This completes the proof of
(\ref{v*}).

In this way we have shown that $U_0$ and $H$ satisfy all the
conditions of Theorem~\ref{T-4}. Therefore, $H$ is the stabilizer
of a functional Menger $\curlywedge$-algebra
$\mathcal{G}_{\curlywedge}$.
 \end{proof}

\begin{minipage} {60mm}
\begin{flushleft}
Dudek~W. A. \\
 Institute of Mathematics and Computer Science \\
 Wroclaw University of Technology \\
 50-370 Wroclaw \\
 Poland \\
 E=mail: dudek@im.pwr.wroc.pl
\end{flushleft}
\end{minipage}
\hfill
\begin{minipage} {60mm}
\begin{flushleft}
 Trokhimenko~V. S. \\
 Department of Mathematics \\
 Pedagogical University \\
 21100 Vinnitsa \\
 Ukraine \\
 E-mail: vtrokhim@sovamua.com
 \end{flushleft}
 \end{minipage}

\end{document}